# Stable and Scalable Probabilistic Numerical Solvers for Stiff and High-Dimensional ODEs

**Nathanael Bosch**
EPFL, Lausanne, Switzerland
nathanael.bosch@epfl.ch

## Abstract

Filtering-based probabilistic numerical solvers for ordinary differential equations (ODEs) have been established as a flexible and efficient simulation framework with built-in numerical uncertainty quantification. However, problems that are both stiff and high-dimensional remain a challenge, as current methods are either stable and have cubic cost in the ODE dimension, or scale linearly at the expense of stability. In this paper, we close this gap and develop probabilistic ODE solvers that are both stable and scalable. We propose two complementary strategies. First, we develop a matrix-free update step that uses Jacobian-vector products, iterative linear solvers, and stochastic covariance estimation to enable linear scaling, all while retaining stability. Second, we propose iterative re-linearization to further improve stability without sacrificing scalability, turning probabilistic ODE solvers into fully implicit methods. We evaluate the proposed approaches on a range of stiff and high-dimensional problems and demonstrate improved stability and scalability over established probabilistic solvers.

## 1 Introduction

Partial differential equations (PDEs) are a fundamental tool used throughout science and engineering to model spatio-temporal dynamical systems. One common approach for simulating such systems is the *method of lines* [22, 29]: by discretizing the PDE spatially, it turns into a system of ordinary differential equations (ODEs), which is then integrated in time with a numerical ODE solver [13]. For many PDEs, however, fine discretization makes the resulting ODE not only high-dimensional but also *stiff*, making it particularly hard to integrate [14]. In this paper, we address this challenge within the framework of *probabilistic numerics* [15, 16, 25] and develop probabilistic numerical ODE solvers for the stiff and high-dimensional regime.

Probabilistic numerical ODE solvers based on Bayesian filtering and smoothing, also known as *ODE filters*, have been established as an efficient and flexible class of probabilistic solvers for ODEs with built-in uncertainty quantification [18, 31, 35]. They have known polynomial convergence rates [18, 36], are applicable to a wide range of differential-equation problems [7], support adaptive step-size selection [5], can be parallelized in time [4], and integrate into larger probabilistic inference tasks [17, 30, 34]. However, for stiff and high-dimensional ODEs current ODE filters face a stability-cost trade-off. Existing A- and L-stable methods have cubic $\mathcal{O}(d^3)$ runtime in the ODE dimension $d$ [6, 35], and ODE filters with linear $\mathcal{O}(d)$ runtime achieve it by approximating the ODE Jacobian and thereby giving up stability [19]. This trade-off also exists within non-probabilistic ODE solvers, where the less-stable *explicit* methods have $\mathcal{O}(d)$ cost whereas stable *implicit* methods scale $\mathcal{O}(d^3)$ [14]. But while in classical numerics this has been addressed with tools like Krylov methods and sparse or structured Jacobians [14, 28], these approaches have not yet been extended to ODE filtering.

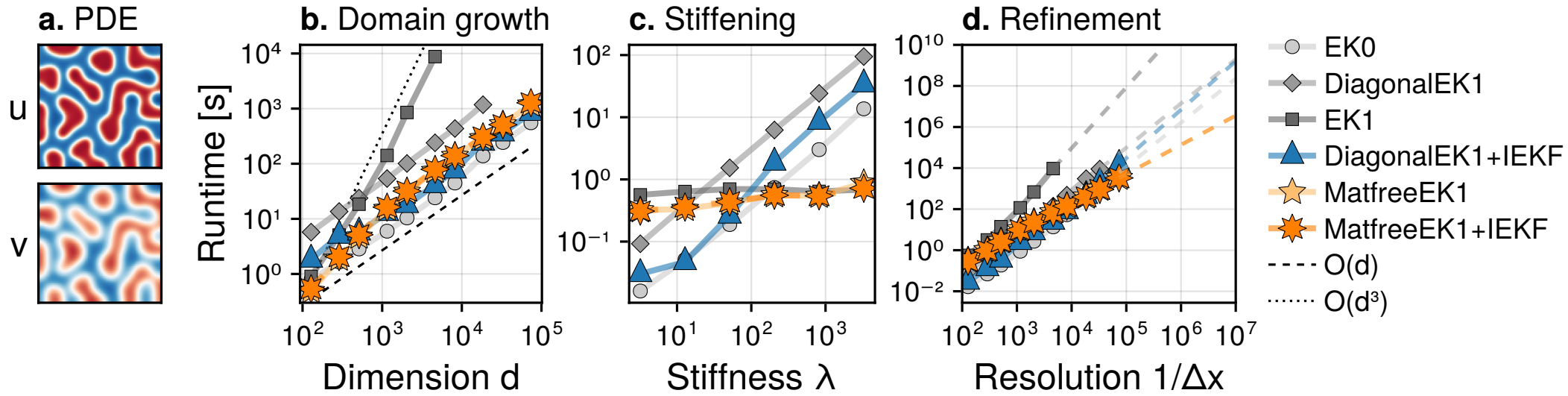

Figure 1: *`MatfreeEK1` is both $\mathcal{O}(d)$ and stable, making it effective for stiff high-dimensional ODEs.* (a) FitzHugh–Nagumo PDE with dimension $d = 125\,000$, solved with `MatfreeEK1`. (b)–(d) Runtime as a function of $d$ at constant stiffness (b), stiffness $\lambda$ at fixed $d$ (c), and PDE resolution $1/\Delta x$ with both $d$ and stiffness increasing (d). Established ODE filters are *either* scalable ($\mathcal{O}(d)$ in (b)) *or* stable (flat in (c)). `MatfreeEK1` achieves both and therefore scales best in (d). See Section A.3.6 for details.

In this paper, we combine tools from classical stable ODE solvers with the covariance approximations of scalable ODE filters to develop probabilistic ODE solvers that are both stable and scalable. We propose two complementary strategies. First, by using Jacobian-vector products, iterative solvers, structured priors as in [19], and stochastic covariance estimation, we formulate a matrix-free update step with linear $\mathcal{O}(d)$ cost, without having to rely on ODE Jacobian approximations. Second, by iteratively re-linearizing as done in *iterated extended Kalman filtering* [12, 33], we develop *fully implicit* ODE filters to further improve stability without affecting the per-step scaling.

**Contributions.** Our main contribution is `MatfreeEK1`, a probabilistic ODE filter that scales linearly in the ODE dimension while retaining the full Jacobian (Section 3.1). We further introduce iterated extended Kalman filtering for ODE filters and to obtain linearly-scaling fully implicit solvers (Section 3.2). We evaluate the proposed methods on a range of stiff and high-dimensional problems and demonstrate their computational scaling, stability under increasing stiffness, and performance at both fixed and adaptive step-size selection (Section 4)

## 2 Probabilistic numerical ODE solvers

Consider an initial value problem with some general non-linear ODE of the form

$$\dot{y}(t) = f(y(t), t), \qquad t \in [0, T], \tag{1}$$

with initial value $y(0) = y_0 \in \mathbb{R}^d$, vector field $f : \mathbb{R}^d \times \mathbb{R} \to \mathbb{R}^d$, and time span $[0, T] \subset \mathbb{R}$. In probabilistic numerics, our goal is to compute a posterior over $y(t)$ such that $y(t)$ satisfies the ODE on a discrete set of points $\mathbb{T} = \{t_n\}_{n=1}^N \subset [0, T]$, that is

$$p\big(y(t) \,\big|\, y(0) = y_0, \{\dot{y}(t_n) = f(y(t_n), t_n)\}_{n=1}^N\big). \tag{2}$$

We call this quantity, and approximations thereof, a *probabilistic numerical ODE solution*, and we call methods that compute such approximations *probabilistic numerical ODE solvers*.

In the following, we provide a brief overview over filtering-based probabilistic numerical ODE solvers, or *ODE filters* [18, 31, 35]. We introduce their underlying prior, observation model, and approximate inference scheme, and then review existing ODE filter variants and their different stability and complexity properties. In Section 3, we then extend this framework to develop probabilistic ODE solvers that are both stable and scalable.

### 2.1 Gauss–Markov prior

*A priori*, we model the ODE solution $y(t)$ as a $q$-times integrated Wiener process (IWP), defined as the output of a linear time-invariant stochastic differential equation (SDE)

$$x(0) \sim N(x(0); \mu_0, \Sigma_0), \qquad \mathrm{d}x = Ax(t)\,\mathrm{d}t + \kappa B\,\mathrm{d}w(t). \tag{3}$$

Here, $x(t) \in \mathbb{R}^D$ is the *state* of the process at time $t \in \mathbb{R}_{\geq 0}$, $A \in \mathbb{R}^{D \times D}$ and $B \in \mathbb{R}^{D \times d}$ are the *drift* and *dispersion* matrices that define the dynamics of the IWP, $w(t) : \mathbb{R}_+ \to \mathbb{R}^d$ is a standard Wiener process, and $\kappa > 0$ is a diffusion parameter. The ODE solution $y(t)$ and its derivative $\dot{y}(t)$, can then be obtained from $x(t)$ via multiplication with suitable *selection matrices* $E_0, E_1 \in \mathbb{R}^{d \times D}$, as $y(t) = E_0 x(t)$ and $\dot{y}(t) = E_1 x(t)$. The IWP-specific values for $A, B, E_0, E_1$ are given in Section A.1.1. To satisfy the initial value, we further assume that $E_0 \mu_0 = y_0$ and $E_0 \Sigma_0 E_0^T = 0$.

In discrete time, the IWP satisfies linear Gaussian transition densities

$$x(t + \Delta t) \mid x(t) \sim N(x(t + \Delta t); \Phi(\Delta t)x(t), Q(\Delta t)), \tag{4}$$

where the matrices $\Phi(\Delta t)$ and $Q(\Delta t)$ can be efficiently computed from $A$ and $B$ [32]. For the IWP, these matrices are known in closed form [18]; formulas are provided in Section A.1.1. Note that other Gauss–Markov processes than the IWP are also admissible and have been used in practice, such as the Matérn process [36] or the integrated Ornstein–Uhlenbeck process [6], but the IWP is the de-facto standard prior choice in ODE filtering due to its favorable computational properties. For a more thorough exposition on Gauss–Markov processes, we refer to [32].

### 2.2 ODE observation model

To connect the prior to the ODE of interest, ODE filters introduce an *observation model* of the form

$$z_n \mid x(t_n) \sim \delta(h_n(x(t_n))), \qquad \text{where} \quad h_n(x) := E_1 x - f(E_0 x, t_n). \tag{5}$$

Under this model, $x(t)$ satisfies the ODE at time $t_n$ if and only if the observation $z_n$ is exactly zero. Conversely, by "observing" $z_n = 0$, we condition $x(t)$ to satisfy the ODE at time $t_n$. This part of the model is also referred to as the *information operator* [9, 36].

### 2.3 Approximate Gaussian inference with extended Kalman filtering

Together, the prior, observation model, and data form a probabilistic inference task, consisting of a non-linear Gaussian state-space model of the form

$$x(0) \sim N(x(0); \mu_0, \Sigma_0), \tag{6a}$$
$$x(t_n) \mid x(t_{n-1}) \sim N(x(t_n); \Phi(\Delta t_n)x(t_{n-1}), Q(\Delta t_n)), \tag{6b}$$
$$z_n \mid x(t_n) \sim \delta(E_1 x(t_n) - f(E_0 x(t_n), t_n)), \tag{6c}$$

and observations $z_n = 0$ for $n = 1, ..., N$. The goal is now to compute a posterior distribution over the state $x(t)$ given the observations $z_{1:N}$, which can then be projected back to the ODE solution space to obtain the desired probabilistic numerical ODE solution. The de-facto standard method used in ODE filtering is the *extended Kalman filter* (EKF):

Starting from the initial distribution $x(0) \sim N(x(0); \mu_0, \Sigma_0)$, the EKF iterates through $t_0, ..., t_N$ and at each step $t_n$ computes Gaussian approximations

$$p(x(t_n) \mid z_{1:n-1}) \approx N(x(t_n); \mu_n^-, \Sigma_n^-), \tag{7a}$$
$$p(x(t_n) \mid z_{1:n}) \approx N(x(t_n); \mu_n, \Sigma_n), \tag{7b}$$

by performing a *prediction*, *linearization*, and *update* step as follows:

- **Prediction:** Extrapolate from $t_{n-1}$ to $t_n$ by marginalizing $x(t_n) \sim N(\mu_n^-, \Sigma_n^-)$ through the linear Gaussian transition model, giving

$$\mu_n^- = \Phi(\Delta t)\mu_{n-1}, \qquad \Sigma_n^- = \Phi(\Delta t)\Sigma_{n-1}\Phi(\Delta t)^T + Q(\Delta t). \tag{8}$$

- **Linearization:** Linearize the observation model $h_n$ by performing a first-order Taylor expansion around the mean estimate $\mu_n^-$, that is, $h_n(x) \approx \hat{h}_n(x) := H_n x + b_n$, where $H_n \in \mathbb{R}^{d \times D}$ is the Jacobian of $h_n$ evaluated at $\mu_n^-$, and $b_n = h_n(\mu_n^-) - H_n \mu_n^- \in \mathbb{R}^d$ is the corresponding offset. For the ODE-filter-specific observation model as defined in (5), we have

$$H_n = E_1 - F_y \cdot E_0, \qquad b_n = F_y E_0 \mu_n^- - f(E_0 \mu_n^-, t_n), \tag{9}$$

where $F_y$ is the Jacobian of the ODE vector field $f$ evaluated at $E_0 \mu_n^-$.

- **Update:** Perform an exact Gaussian update on the linearized observation model and compute

$$\mu_n = \mu_n^- + K_n\big(0 - \bar{h}_n(\mu_n^-)\big), \qquad \Sigma_n = (I - K_n H_n)\Sigma_n^-(I - K_n H_n)^T, \tag{10}$$

with $K_n := \Sigma_n^- H_n^T \left(H_n \Sigma_n^- H_n^T\right)^{-1}$.

After the EKF has run through all time points $t_{1:N}$ one can optionally perform a *smoothing* backward pass to refine the solution; see Section A.1.2 for details. We then project the state posteriors $N(x(t_n); \mu_n, \Sigma_n)$ back to the ODE solution space via multiplication with the selection matrix $E_0$ and obtain the desired probabilistic numerical ODE solution:

$$p\Big(y(t_n) \mid y(0) = y_0, \{\dot{y}(t_n) = f(y(t_n), t_n)\}_{n=1}^N\Big) \approx N(y(t_n); E_0\mu_n, E_0\Sigma_n E_0^T). \tag{11}$$

This is precisely the filtering-based probabilistic numerical ODE solver as introduced by [35].

### 2.4 ODE filter variants

This framework forms the basis of a range of ODE filtering algorithms, with different stability and complexity properties. They primarily differ in how the ODE is linearized within the solver.

- **`EK1`:** Uses the exact ODE Jacobian, $F_y = \partial_y f(E_0\mu_n^-, t_n)$, as described in the procedure above. The resulting method is *A-stable* [35] and can be considered a *linearly implicit* or *semi-implicit* ODE solver [3, 14], making it well-suited for stiff problems. Due to the required matrix multiplications and inversions in the predict and update steps, a single `EK1` step costs $\mathcal{O}(d^3)$, and requires $O(d^2)$ memory to store the dense covariance and Jacobian matrices.
- **`EK0`:** Sets $F_y = 0$, resulting in an *explicit* ODE filter without any stability guarantees [31]. When the chosen prior has Kronecker structure (as is the case for the IWP), the `EK0` prediction and update preserve this structure and the method achieves $\mathcal{O}(d)$ scaling in time and space [19].
- **`DiagonalEK1`:** Approximates the Jacobian by its diagonal, $F_y \approx \operatorname{diag}\big(\partial_y f\big)$ [19]. The `DiagonalEK1` preserves block-diagonal structure in the prior (present in the IWP) and thereby also achieves $\mathcal{O}(d)$ scaling. Compared to `EK0`, the partial Jacobian empirically improves stability [19].
- **`ExpEK`:** Takes a different approach and modifies the prior rather than the update, in such a way that the linearized ODE is solved *exactly* during the prediction step. This makes the method *L-stable* [6]. However, `ExpEK` does not preserve any covariance structure and thus falls back to $\mathcal{O}(d^3)$ time and $\mathcal{O}(d^2)$ memory cost.

Current ODE filters are therefore either stable at $\mathcal{O}(d^3)$ time and $\mathcal{O}(d^2)$ memory cost (`EK1`, `ExpEK`), or are $\mathcal{O}(d)$ at the expense of stability (`EK0`, `DiagonalEK1`).

## 3 Improving the stability of scalable ODE filters

Here we develop probabilistic ODE solvers that are both stable on stiff problems and scalable to high-dimensional systems, relying on structured priors but without Jacobian approximations. We pursue two complementary strategies. First, we propose a *matrix-free* EKF update step that linearizes exactly without instantiating the Jacobian, at $\mathcal{O}(d)$ cost per step (Section 3.1). Second, we propose *iterative re-linearization* to further improve stability without sacrificing scalability, turning the solver into a fully implicit method (Section 3.2).

### 3.1 Efficient EKF updates without Jacobian instantiations

Recall the EKF update equations from Section 2.3, given by

$$\mu_n = \mu_n^- - K_n(H_n\mu_n^- + b_n), \tag{12a}$$

$$\Sigma_n = (I - K_n H_n)\Sigma_n^-(I - K_n H_n)^T, \tag{12b}$$

with $H_n$, $b_n$, and $K_n$ given by

$$H_n = E_1 - F_y E_0, \qquad b_n = F_y E_0 \mu_n^- - f(E_0\mu_n^-, t_n), \qquad K_n = \Sigma_n^- H_n^T \left(H_n \Sigma_n^- H_n^T\right)^{-1}. \tag{13}$$

In general, performing these computations comes with $\mathcal{O}(d^3)$ computational cost due to the required matrix inversion and matrix-multiplications. But as discussed in Section 2.4, the `DiagonalEK1` achieves $\mathcal{O}(d)$ cost per step by combining a block-diagonal prior covariance with a diagonal Jacobian approximation $F_y \approx \operatorname{diag}(\partial_y f)$: the Jacobian approximation preserves block-diagonality, and the block-diagonality in turn enables efficient $\mathcal{O}(d)$ matrix-products and matrix inversions. In this section, we develop an update that uses the exact Jacobian and instead directly imposes block-diagonality on the posterior covariance $\Sigma_n$ to achieve $\mathcal{O}(d)$ cost per step.

### 3.1.1 Matrix-free EKF mean updates

Tracing the mean update $\mu_n^- \to \mu_n$ through the equations above, we observe that the (potentially) costly computations are (i) matrix-vector multiplications with $E_0, E_1$, $F_y$, $\Sigma_n^-$ and transposes, (ii) inverting $\left(H_n \Sigma_n^- H_n^T\right)$ and (iii) evaluating the vector-field $f$. The latter is unavoidable in any numerical ODE solver; we denote its cost by $C_{\mathrm{f}}$. Regarding (i) and (ii), we can leverage the following:

(i) *Matrix-vector products:* The equations contain three types of matrix-vector products. First, multiplication with $E_0, E_1$ can be implemented efficiently by indexing operations since $E_0, E_1$ are selection matrices, and is thus negligible [19]; see also Section A.1.1. Second, multiplication with $\Sigma_n^-$ has cost $\mathcal{O}(d^2)$ in general, but as we assume block-diagonal structure its cost reduces to $\mathcal{O}(d)$ [19]. Third, multiplication with $F_y$ and $F_y^T$ can be done efficiently without instantiating the Jacobian via automatic differentiation, at cost comparable to evaluating the vector field; these are also known as *Jacobian-vector products* (JVPs) and *vector-Jacobian products* (VJPs).

(ii) *Matrix inversion:* Instead of inverting the matrix $\left(H_n \Sigma_n^- H_n^T\right)$ directly, we can solve the corresponding linear system using an iterative solver such as conjugate gradients (CG). This has cost $\mathcal{O}(k C_{\mathrm{mv}})$, where $k$ is the number of solver iterations and $C_{\mathrm{mv}}$ is the cost of performing a matrix-vector product with the system matrix. As discussed in (i), we have $C_{\mathrm{mv}} = \mathcal{O}(d + C_{\mathrm{f}})$.

Putting everything together, JVPs/VJPs, iterative linear solvers, and block-diagonal covariances enable *exact* EKF mean updates at cost $\mathcal{O}(k \cdot (d + C_{\mathrm{f}}))$, without requiring approximate Jacobians.

### 3.1.2 Stochastically-estimated block-diagonal covariance updates

To preserve the computational benefits of the block-diagonal structure, we develop an update step that approximates the "correct" posterior covariance $\Sigma_n$ by a stochastically-estimated block-diagonal matrix in the same factorization as the prior covariance $\Sigma_n^-$, again in a matrix-free manner.

We first rewrite the EKF covariance update in *square-root form* [20] as

$$\sqrt{\Sigma_n} = (I - K_n H_n) \sqrt{\Sigma_n^-}, \tag{14}$$

where $\sqrt{M}$ denotes a *matrix square-root* of $M$, defined such that $\sqrt{M}\left(\sqrt{M}\right)^T = M$.

Next, we generate $m$ independent standard normal samples $\eta_1, ..., \eta_m \sim N(0, I_D)$, and then linearly transforming these samples with $\xi_k = \sqrt{\Sigma_n} \eta_k$, that is,

$$\xi_k = (I - K_n H_n) \sqrt{\Sigma_n^-} \eta_k, \qquad \text{for } k = 1, ..., m. \tag{15}$$

This can be computed efficiently at cost $\mathcal{O}(k \cdot (d + C_{\mathrm{f}}))$ via matrix-vector products with $K_n$, $H_n$, and $\sqrt{\Sigma_n^-}$, ; see Section 3.1.1. We obtain $m$ independent samples with distribution $\xi_k \sim N(0, \Sigma_n)$.

Finally, we compute a block-diagonal approximation of $\Sigma_n$ from these samples: Let $\imath$ be the index set of the $i$-th block of $\Sigma_n$ such that $[\Sigma_n]_{\imath,\imath} \in \mathbb{R}^{(q+1)\times(q+1)}$. We then approximate each block with

$$[\Sigma_n]_{\imath,\imath} \approx \frac{1}{m} \sum_{k=1}^{m} (\xi_k)_\imath (\xi_k)_\imath^T. \tag{16}$$

This directly follows from the identity $\mathbb{E}[\xi \xi^T] = \Sigma$ for $\xi \sim N(0, \Sigma)$. Since $(\xi_k)_\imath \in \mathbb{R}^{q+1}$, the cost of computing each of the $d$ blocks is $\mathcal{O}(m(q+1)^2)$. Note that it is also possible to directly compute a square-root factor $\sqrt{[\Sigma_n]_{\imath,\imath}}$ from the samples at the same cost; see Section A.2 for details.

We obtain a block-diagonal covariance update at cost $\mathcal{O}\big(m(q+1)^2 + k \cdot (d + C_{\mathrm{f}})\big)$. Assuming that the number of linear solver iterations $k$ and the number of samples $m$ are independent of $d$, the update thus scales linearly in the ODE dimension $d$. The block-diagonal structure also keeps the memory cost at $\mathcal{O}(d)$.

Together with the matrix-free mean update of Section 3.1.1, this completes the matrix-free $\mathcal{O}(d)$ EKF update. We refer to the resulting ODE filter variant as `MatfreeEK1`.

### 3.2 Iterative re-linearization

To further improve stability, we propose to iteratively refine the linearization point at each update step, using what is known as the *iterated extended Kalman filter* (IEKF) in Bayesian filtering and smoothing [12, 33]. The IEKF can be understood as a Gauss–Newton algorithm that computes the maximum-a-posteriori (MAP) estimate

$$\mu_n = \operatorname*{argmax}_x p(x \mid z_{1:n}) \tag{17}$$

at each time point [1]. Since this describes an implicit equation involving $f$, a probabilistic numerical ODE solver that uses the IEKF is *fully implicit* in the classical sense [14].

The IEKF operates as follows. Starting from the prediction mean $\xi = \mu_n^-$, the IEKF iteratively performs *linearization* and *update* steps until a *stopping criterion* is met:

- **Linearization:** Approximate the non-linear observation model by Taylor-expanding it around the current linearization point $\xi$; this can be done with the exact Jacobian, a diagonal Jacobian approximation, or matrix-free as in Section 3.1.
- **Update:** Using the linearized observation model, compute the posterior mean $\mu_n^{(k)}$; this can again be done with an exact EKF update, or with the update proposed in Section 3.1.1.
- **Stopping criterion:** If $\|\mu_n^{(k)} - \xi\| < \varepsilon$ for some chosen $\varepsilon > 0$, set $\mu_n \leftarrow \mu_n^{(k)}$ and stop iterating. If not, update the linearization point to the latest updated mean $\xi \leftarrow \mu_n^{(k)}$ and continue.

After the mean has converged, compute the corresponding covariance estimate $\Sigma_n$ under the observation model linearized at the converged mean, using any of the update variants above.

The IEKF does not affect the per-step scaling of the underlying solver, as long as the iteration count is constant with respect to the ODE dimension $d$. In principle, it could be combined with any of the variants described in Section 2.4. We focus on combinations with the `DiagonalEK1` and the proposed `MatfreeEK1` to obtain ODE filters that are both stable and scalable.

### 3.3 Implementation details

We solve the linear system within the mean update using conjugate gradients (CG), and solve the $m$ linear systems of the covariance sample computation jointly via block-GMRES [28]. For the stochastic block-diagonal estimation, we default to $m = 2(q+1)$ samples. All filtering and smoothing operations are implemented in square-root form as described in [20]. We additionally perform posterior calibration and step-size adaptation as proposed by [5]. Default tolerances and further implementation details are provided in Section A.2.

## 4 Experiments

We evaluate the proposed `MatfreeEK1`, `MatfreeEK1+IEKF`, and `DiagonalEK1+IEKF` against the existing probabilistic ODE solvers from Section 2.4 on a range of stiff and high-dimensional problems, and analyze computational scaling, stability under increasing stiffness, and performance under fixed and adaptive step-size selection.

All methods are implemented in the Julia programming language [2]. JVPs and VJPs were computed with ForwardDiff.jl [27] and Enzyme.jl [24] respectively, through DifferentiationInterface.jl [10]. Iterative solvers were provided by Krylov.jl [23]. Reference solutions were computed with

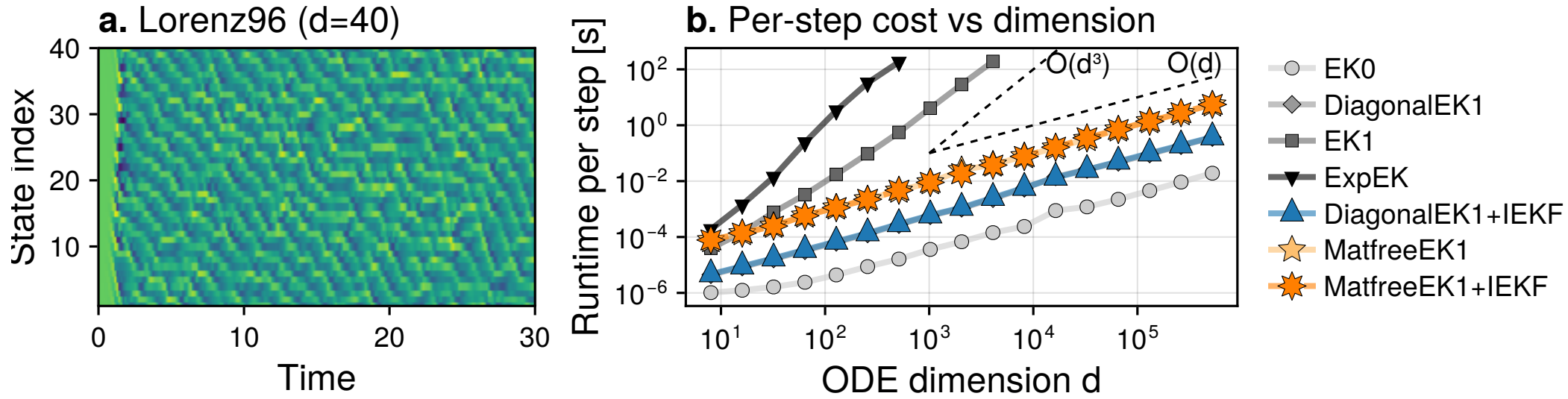

Figure 2: *`MatfreeEK1` scales linearly in the ODE dimension.* Wall-clock time of a single solver step on the Lorenz96 system as a function of dimension $d$. `MatfreeEK1` matches the linear scaling of `EK0` and `DiagonalEK1`; `EK1` and `ExpEK` scale cubically. `IEKF` does not affect scaling behavior.

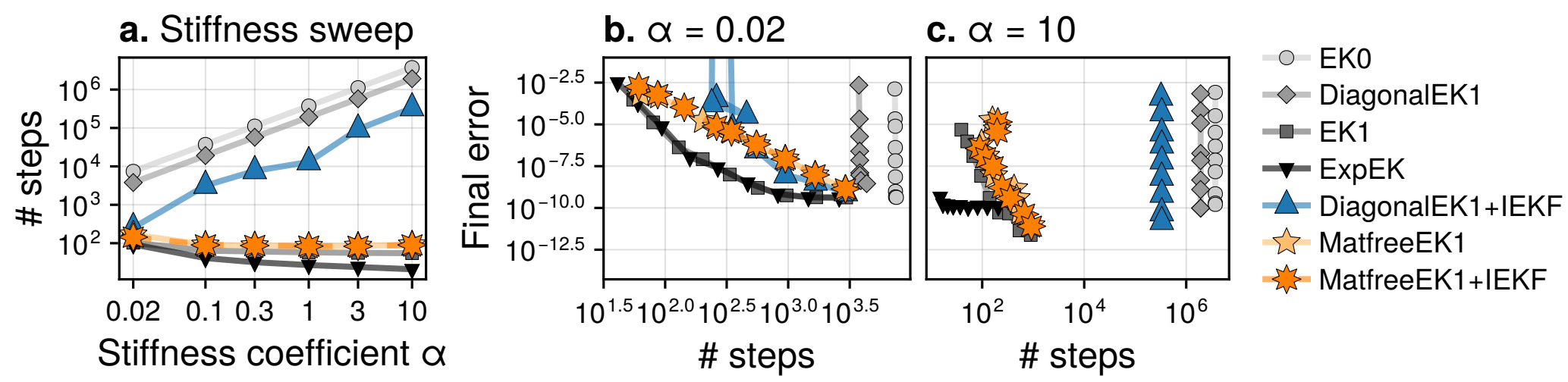

Figure 3: *`MatfreeEK1` is stable under increasing stiffness.* Number of steps as a function of the stiffness parameter $\alpha$ of the 1D Brusselator (a), and work-precision diagrams for the least and most stiff settings (b, c). `EK0` and `DiagonalEK1` require more steps as stiffness increases, whereas `MatfreeEK1`, `EK1`, and `ExpEK` keep their step count essentially constant. `IEKF` does not affect `MatfreeEK1` but reduces the step count of `DiagonalEK1`.

DifferentialEquations.jl [26]. All experiments ran on an *AMD Ryzen Threadripper PRO 7985WX* CPU. Code for the implementation and experiments will be made publicly available on GitHub.[1]

### 4.1 Computational scaling on Lorenz96

We first evaluate the per-step computational cost of the proposed methods as the ODE dimension increases. Following the setup of [19], we consider the Lorenz96 system [21], which is a chaotic dynamical system whose dimension $d$ can be freely chosen without affecting stiffness. For a range of dimensions $d$, we evaluate the runtime of a single step of the proposed methods `MatfreeEK1`, `MatfreeEK1+IEKF`, and `DiagonalEK1+IEKF` and compare them to `EK0`, `EK1`, `DiagonalEK1`, and `ExpEK`. The full experimental setup is detailed in Section A.3.

The results are shown in Figure 2. `MatfreeEK1` matches the linear scaling of `EK0` and `DiagonalEK1`, while `EK1` and `ExpEK` scale cubically. We also observe that the `IEKF` does not change the scaling behavior of the method it is applied to. This confirms the previously discussed computational complexity of the proposed methods and establishes `MatfreeEK1` as a scalable ODE filter.

### 4.2 Stability under increasing stiffness on the 1D Brusselator

Next, we evaluate the stability of the proposed methods under increasing stiffness. We consider the 1D Brusselator [14], which is a reaction-diffusion system with stiffness parameter $\alpha$, and solve it across a range of $\alpha$ values. For each value, we report the number of steps required to solve the problem with adaptive step-size selection given a fixed target tolerance. As before, we evaluate the proposed methods and compare them to `EK0`, `EK1`, `DiagonalEK1`, and `ExpEK`. The full experimental setup is detailed in Section A.3.

[1] https://github.com/nathanaelbosch/stable-scalable-probnum-ode-solvers

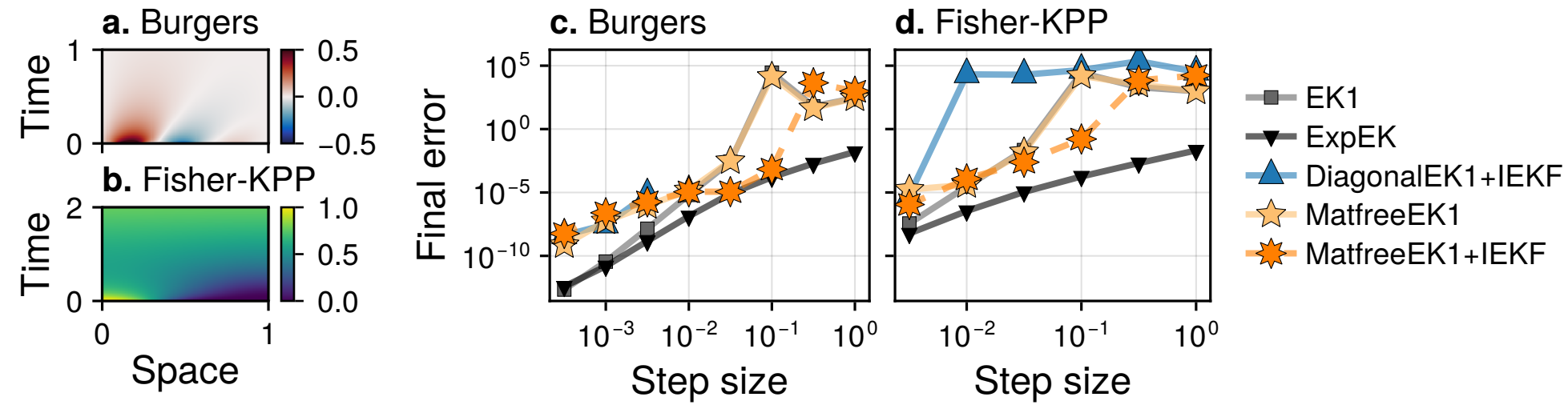

Figure 4: *Fixed-step performance of probabilistic ODE solvers on stiff PDEs.* ODE solutions of Burgers equation and Fisher-KPP (a, b), and final error vs. step size for both problems (c, d). `MatfreeEK1` shows the similar stability as `EK1`. `IEKF` improves stability of `MatfreeEK1` at large $\Delta t$.

Figure 3 shows the results. We observe that as $\alpha$ grows and stiffness increases, `EK0` and `DiagonalEK1` require more steps to maintain the same accuracy while `MatfreeEK1`, `EK1`, and `ExpEK` keep their step counts essentially constant. We also see that `IEKF` reduces the step count of the `DiagonalEK1`, but does not change its qualitative behaviour with respect to stiffness. The `IEKF` does not affect `MatfreeEK1`, suggesting that its exact linearization is already sufficient for stability.

For the least and most stiff settings we also report work-precision diagrams in Figure 3. We observe a widening gap between the stable and unstable methods as stiffness increases, as well as vertical curves for `EK0`, `DiagonalEK1`, and `DiagonalEK1+IEKF`. This indicates that small steps need to be taken independently of the specified tolerance and is known as *stability-limited*, a common behaviour of explicit methods on stiff problems [14]. In contrast, `MatfreeEK1` and `EK1` maintain an accuracy-cost trade-off across most of the tolerance range. `ExpEK` on the other hand becomes almost horizontal at $\alpha = 10$, as it is able to resolve the dominating linear dynamics in closed-form [6]. This experiment establishes the proposed `MatfreeEK1` as a stable solver for stiff problems.

### 4.3 Stability at fixed step sizes on stiff PDEs

To isolate stability from step-size adaptation, we next evaluate the proposed methods at fixed step sizes. Following the setup of [6] we consider two stiff semi-linear PDEs, Burgers and Fisher-KPP, and compare the proposed `MatfreeEK1`, `MatfreeEK1+IEKF`, and `DiagonalEK1+IEKF` against `EK1` and `ExpEK` across a range of fixed step sizes. `EK0` and `DiagonalEK1` are omitted as they fail to produce meaningful results at any considered step size. The full experimental setup is detailed in Section A.3.

Figure 4 shows the results. `MatfreeEK1` agrees closely with `EK1` at larger step sizes $\Delta t \geq 10^{-2}$, but deviates for smaller steps. `ExpEK` is the most stable across the $\Delta t$ range, but has $\mathcal{O}(d^3)$ cost as established in [6]. We also observe a stability benefit of the `IEKF` for `MatfreeEK1` at large step sizes, indicating that in the absence of adaptive step-size control, re-linearization can help maintain stability. The `IEKF` also improves the stability of `DiagonalEK1`, which otherwise fails for all considered step sizes. This experiment further confirms the stability of `MatfreeEK1` and shows that, at fixed step sizes, iterative re-linearization improves stability.

### 4.4 Adaptive solving on stiff PDEs

Finally, we evaluate the proposed methods with adaptive step-size selection on two high-dimensional discretized stiff PDEs. We consider 1D Burgers equation together with a 2D Fisher-KPP model, resulting in ODEs of dimension $d = 500$ and $d = 4096$ respectively. We evaluate the proposed `MatfreeEK1`, `MatfreeEK1+IEKF`, and `DiagonalEK1+IEKF` against `EK1` and `DiagonalEK1` on a range of tolerances, skipping `ExpEK` due to high computational cost. The full setup is detailed in Section A.3.

Figure 5 shows the results. On both equations, `MatfreeEK1` shows an accuracy-cost trade-off across the whole tolerance range, whereas `DiagonalEK1` and `DiagonalEK1+IEKF` show stability-limited behaviour with vertical curves. `MatfreeEK1` also matches the convergence slope of `EK1` on Burgers, although with a slight offset; since both methods share the same mean update, we suspect that this relates to the covariance approximation within `MatfreeEK1`. We again observe no effect of `IEKF` on `MatfreeEK1`, but a noticeable reduction in runtime for `DiagonalEK1`, turning it a competitive option at low tolerances. Note that on Fisher-KPP, the $O(d^3)$ cost of `EK1` lead to ~1h45m runtime at the loosest tolerance and was thus omitted from the plot. `MatfreeEK1` on the other hand was able to solve the problem across the whole tolerance range, with a runtime of 27min at the tightest tolerance.

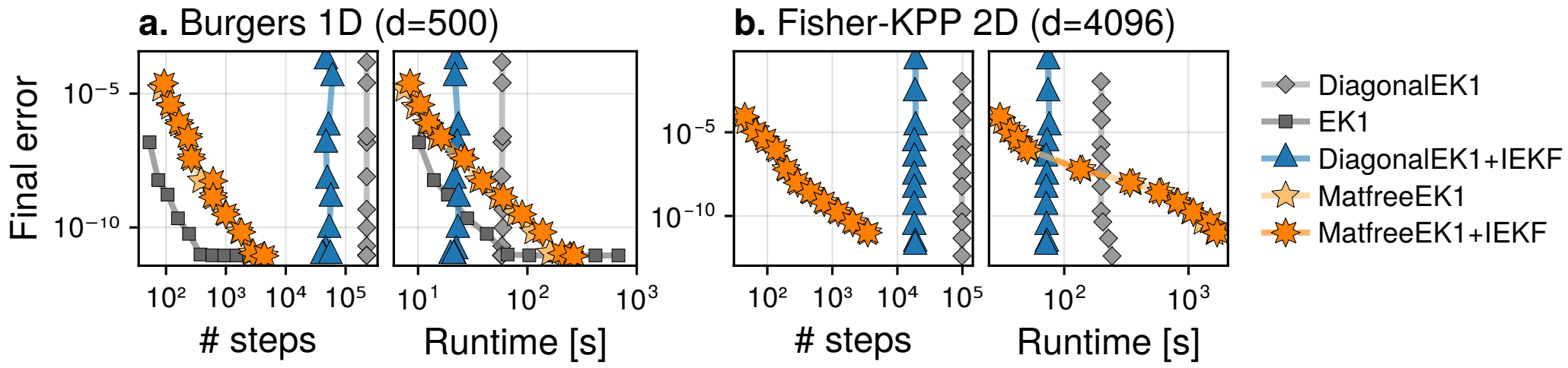

Figure 5: *Benchmarking ODE filters on large stiff PDEs.* Final error vs. step count and wall-clock time on Burgers ($d = 500$) and Fisher-KPP ($d = 4096$). `MatfreeEK1` shows an accuracy-cost trade-off across the whole tolerance range, while `DiagonalEK1` (with and without `IEKF`) is stability-limited on both problems. `EK1` is omitted on Fisher-KPP where its $\mathcal{O}(d^3)$ cost is infeasible.

Altogether, these experiments establish `MatfreeEK1` as a stable and scalable probabilistic ODE solver and the `IEKF` as a useful accelerator for the existing `DiagonalEK1`.

## 5 Limitations

**Utility of iterative re-linearization.** Section 4.3 shows that `MatfreeEK1` benefits from `IEKF` iterations at fixed step sizes, but this improvement did not extend to the adaptive setting in Section 4.4. The `DiagonalEK1` however did benefit from `IEKF` in both fixed and adaptive settings. This suggests that step-size adaptation provides much of the `IEKF` benefit when the Jacobian is accurate, but a more thorough investigation could help clarify when re-linearization is most useful. A comparison to the non-scalable `EK1+IEKF` combination was beyond the scope of this work but could help further clarify the utility of iterative re-linearization.

**Stability of `MatfreeEK1`.** While the `EK1` and `ExpEK` are known to be A- and L-stable respectively, the proposed `MatfreeEK1` does not yet have an analogous theoretical guarantee as the A-stability proof of [35] does not directly transfer due to the posterior covariance approximation. The effect of the number of covariance estimator samples was also not studied systematically. We leave a more thorough theoretical analysis of stability properties of the `MatfreeEK1` to future work.

**Posterior calibration.** While we have applied calibration to enable adaptive step-size selection [5], our experiments report only point-estimate metrics and do not directly evaluate the uncertainty of the resulting probabilistic ODE solution. Since evaluating whether a probabilistic numerical method is well-calibrated is itself a non-trivial research question [8], we leave a thorough investigation of posterior uncertainties and their calibration to future work.

## 6 Conclusion

Probabilistic numerical ODE solvers have so far had to choose between stability for stiff systems and scalability on high-dimensional problems. We proposed two complementary strategies to close this gap. `MatfreeEK1` combines matrix-free linear algebra and stochastic covariance estimation to perform a full-Jacobian update without instantiating it, at $\mathcal{O}(d)$ cost per step. Iterative re-linearization via the `IEKF` improves stability without affecting the per-step cost, and turns the resulting ODE filters into fully implicit methods. Experiments on a range of stiff and high-dimensional systems confirm both of these properties and demonstrate the utility of the proposed methods, in particular for solving large stiff discretized PDEs. Together, these contributions draw new connections between probabilistic solvers and the classical numerical analysis toolbox, extend the reach of probabilistic numerical ODE solvers to a wider range of problems, and open up new avenues for future research on stable and scalable probabilistic numerical methods.

# A Appendix

## A.1 Further details on filtering-based probabilistic numerical ODE solvers

### A.1.1 The $q$-times Integrated Wiener Process Prior

The de-facto most-commonly used Gauss–Markov prior for probabilistic ODE solvers is the *$q$-times integrated Wiener process* (IWP), chosen independently for each dimension of the ODE. This corresponds to SDE matrices of the form $A = \breve{A} \otimes I_d$ and $B = \breve{B} \otimes I_d$ with

$$\breve{A} = \begin{bmatrix} 0 & 1 & \cdots & 0 \\ \vdots & \vdots & \ddots & \vdots \\ 0 & 0 & \cdots & 1 \\ 0 & 0 & \cdots & 0 \end{bmatrix} \in \mathbb{R}^{(q+1)\times(q+1)}, \qquad \breve{B} = \begin{bmatrix} 0 \\ \vdots \\ 0 \\ 1 \end{bmatrix} \in \mathbb{R}^{(q+1)}, \tag{18}$$

and selection matrices

$$\begin{aligned} E_0 &= [1 \ \ 0 \ \ \cdots \ \ 0] \otimes I_d, \\ E_1 &= [0 \ \ 1 \ \ \cdots \ \ 0] \otimes I_d, \end{aligned} \tag{19}$$

and resulting state dimension $D = (q+1)d$. Notably, the Kronecker structure of this dimension-independent IWP prior also enables the favorable computational scaling developed by [19]. We therefore also adopt it as the default choice in this work.

Due to the Kronecker structure of $A$ and $B$, the discrete-time transition matrices $\Phi(h)$ and $Q(h)$ also factor as Kronecker products,

$$\Phi(h) = \breve{\Phi}(h) \otimes I_d, \qquad Q(h) = \kappa^2 \breve{Q}(h) \otimes I_d, \tag{20}$$

where the $(q+1) \times (q+1)$ matrices $\breve{\Phi}(h)$ and $\breve{Q}(h)$ are given in closed form by [18]:

$$\begin{aligned} \left[\breve{\Phi}(h)\right]_{ij} &= \mathbb{1}_{i \leq j} \cdot \frac{h^{j-i}}{(j-i)!}, \\ \left[\breve{Q}(h)\right]_{ij} &= \frac{h^{2q+1-i-j}}{(2q+1-i-j)(q-i)!(q-j)!}, \end{aligned} \tag{21}$$

for $i, j = 0, ..., q$.

### A.1.2 Smoothing

The EKF forward pass produces *filtering* posteriors $p(x(t_n) \mid z_{1:n}) \approx N(x(t_n); \mu_n, \Sigma_n)$ that only condition on observations up to time $t_n$. *Smoothing* is a subsequent backward pass that refines these into *smoothing* posteriors $p(x(t_n) \mid z_{1:N}) \approx N\big(x(t_n); \mu_n^S, \Sigma_n^S\big)$, which condition on all observations [33]. Starting from $\mu_N^S = \mu_N$ and $\Sigma_N^S = \Sigma_N$, the Rauch–Tung–Striebel smoother iterates through $n = N-1, ..., 0$ and computes

$$\begin{aligned} G_n &= \Sigma_n \Phi(\Delta t_n)^T \left(\Sigma_{n+1}^{-}\right)^{-1}, \\ \mu_n^S &= \mu_n + G_n\big(\mu_{n+1}^S - \mu_{n+1}^{-}\big), \\ \Sigma_n^S &= \Sigma_n + G_n\big(\Sigma_{n+1}^S - \Sigma_{n+1}^{-}\big)G_n^T, \end{aligned} \tag{22}$$

where $G_n \in \mathbb{R}^{D \times D}$ is the smoother gain.

In this work, we perform smoothing to obtain probabilistic solution estimates over the full time span $[0, T]$. If only the solution at the final time point $t_N$ is of interest, the smoothing pass can be omitted.

## A.2 Further implementation details of the proposed `MatfreeEK1` and `IEKF`

**Iterative linear solvers.** The mean update linear system is solved with conjugate gradients (CG) via Krylov.jl [23]. The CG tolerance inherits the relative tolerance specified to the ODE solver for its step-size adaptation, capped at $10^{-2}$ to prevent premature convergence at loose tolerances; this is common practice in ODE solvers with iterative linear algebra [26]. The $m$ linear systems of the covariance sample computation are solved jointly using block-GMRES, again via Krylov.jl, with the same tolerance as CG. When an analytic diagonal Jacobian is available, we optionally apply a Jacobi preconditioner based on $\mathrm{diag}(F_y)$ to both linear solves.

**Stochastic covariance estimator.** We use $m = 2(q+1)$ independent Gaussian probe vectors by default to construct the stochastic block-diagonal covariance estimator.

**Iterative re-linearization.** The IEKF is capped at 20 iterations per step, with early stopping when the scaled update step falls below $\varepsilon = 10^{-2}$.

**Square-root implementation.** All filtering and smoothing operations are carried out in square-root form to enhance numerical stability, as is common practice in ODE filtering [20]. The prediction step maintains a square-root factor of the covariance via a QR decomposition. For the stochastic block-diagonal covariance estimator, a square-root factor of each block is obtained directly from the samples without explicitly forming the block. Recall from Section 3.1.2 that the samples $\xi_1, ..., \xi_m \sim N(0, \Sigma_n)$ are available, and that the $i$-th block is estimated as

$$[\Sigma_n]_{i,i} \approx \frac{1}{m} \sum_{k=1}^{m} (\xi_k)_i (\xi_k)_i^T. \tag{23}$$

Stacking the sub-vectors into the matrix

$$S_i = \frac{1}{\sqrt{m}} \begin{bmatrix} (\xi_1)_i^T \\ \vdots \\ (\xi_m)_i^T \end{bmatrix} \in \mathbb{R}^{m \times (q+1)}, \tag{24}$$

we have $S_i^T S_i = \frac{1}{m} \sum_k (\xi_k)_i (\xi_k)_i^T \approx [\Sigma_n]_{i,i}$. A thin QR decomposition $S_i = Q_i R_i$ then gives $R_i^T R_i = S_i^T S_i$, so $R_i^T$ is a square-root factor of the estimated block, computed at cost $\mathcal{O}(m(q+1)^2)$ without ever forming the $(q+1) \times (q+1)$ block itself. We refer to [20] for a detailed treatment of square-root implementations in the ODE filter setting.

### A.3 Problem descriptions

We use the following conventions throughout. 1D PDEs are discretized on an equidistant grid of $N$ points on $\Omega = (0,1)$ with spacing $\Delta x = 1/N$; 2D PDEs analogously on $N \times N$ cells of $\Omega = (0,1)^2$. For experiments with adaptive step-size selection, the *WP tolerance grid* is the paired set $\tau_{\text{rel}} \in \{10^1, 10^0, ..., 10^{-9}\}$, $\tau_{\text{abs}} = 10^{-3} \cdot \tau_{\text{rel}}$.

#### A.3.1 Lorenz96

The Lorenz96 system [21] is a chaotic dynamical system on a ring of $N$ nodes:

$$\dot{y}_i = (y_{i+1} - y_{i-2}) y_{i-1} - y_i + F, \quad i = 1, ..., N, \tag{25}$$

with periodic indexing $y_{i+N} = y_i$, forcing $F = 8$ (chaotic regime), and initial value $y_i(0) = F$ for $i > 1$ and $y_1(0) = F + 0.01$. Dimension $d = N$.

**Further details.** For the scaling experiment of Section 4.1, we sweep $N \in \{2^3, ..., 2^{21}\}$ and time a single solver step at fixed step size $\Delta t = 0.01$ with an IWP(3) prior, reporting the minimum over five timed steps after one warmup step.

#### A.3.2 1D Brusselator

The 1D Brusselator [14] is a two-species reaction-diffusion system on $\Omega = (0,1)$:

$$\begin{aligned} \partial_t u_1(x,t) &= \alpha \partial_x^2 u_1 + A + u_1^2 u_2 - (B+1) u_1, \\ \partial_t u_2(x,t) &= \alpha \partial_x^2 u_2 + B u_1 - u_1^2 u_2, \end{aligned} \tag{26}$$

with parameters $A = 1$, $B = 3$, diffusion coefficient $\alpha > 0$ controlling stiffness, Dirichlet boundary conditions $u_1(0,t) = u_1(1,t) = 1$, $u_2(0,t) = u_2(1,t) = 3$, initial conditions

$$u_1(x,0) = 1 + \sin(2\pi x), \qquad u_2(x,0) = 3, \tag{27}$$

and time span $t \in [0, 10]$. Each species is discretized on a grid of $N = 40$ points using the standard three-point Laplacian

$$[L]_{ij} = \frac{1}{(\Delta x)^2} \cdot \begin{cases} -2 & \text{if } i = j \\ 1 & \text{if } i = j \pm 1 \\ 0 & \text{otherwise,} \end{cases} \tag{28}$$

and enforce the Dirichlet conditions by holding the four boundary entries of the discrete state fixed. Dimension $d = 2N = 80$.

**Further details.** The $\alpha$-sweep panel of Section 4.2 varies $\alpha$ logarithmically over $\{10^{-2}, ..., 10^{1}\}$ at fixed tolerances $\tau_{\text{abs}} = 10^{-6}$, $\tau_{\text{rel}} = 10^{-3}$; the work-precision diagrams at $\alpha = 0.02$ and $\alpha = 10$ use the WP tolerance grid. All methods use an IWP(3) prior with adaptive step-size selection.

#### A.3.3 1D Burgers

The 1D Burgers equation is the semi-linear PDE

$$\partial_t u(x,t) = -u(x,t)\partial_x u(x,t) + D\partial_x^2 u(x,t), \qquad x \in \Omega, \qquad t \in [0,1], \tag{29}$$

with diffusion $D = 0.075$, zero-Dirichlet boundary conditions $u(0,t) = u(1,t) = 0$, and initial condition $u(x,0) = \sin(3\pi x)^3(1-x)^{3/2}$. We discretize using the same three-point Dirichlet Laplacian $L$ as in Section A.3.2 together with a conservative-form discretization of $u \cdot \partial_x u = (1/2)\partial_x(u^2)$,

$$[F(y)]_i = \frac{1}{4\Delta x} \cdot \begin{cases} y_2^2 & \text{if } i = 1 \\ -y_{N-1}^2 & \text{if } i = N \\ y_{i+1}^2 - y_{i-1}^2 & \text{otherwise.} \end{cases} \tag{30}$$

This yields the semi-linear ODE $\dot{y}(t) = DLy(t) - F(y(t))$, where the boundary entries of $L$ and $F$ implicitly enforce the zero-Dirichlet conditions via ghost-zero extension.

**Further details.** The fixed-step experiment of Section 4.3 uses $N = 200$ with an IWP(2) prior and step grid $\Delta t \in \{10^{-k/2}\}_{k=0}^{8}$; the adaptive work-precision experiment of Section 4.4 uses $N = 500$ with an IWP(3) prior and the WP tolerance grid.

#### A.3.4 1D Fisher-KPP

The 1D Fisher-KPP equation [11] is the semi-linear reaction-diffusion PDE

$$\partial_t u(x,t) = D\partial_x^2 u(x,t) + u(x,t)(1 - u(x,t)), \qquad x \in \Omega, \qquad t \in [0,2], \tag{31}$$

with diffusion $D = 0.25$, zero-Neumann boundary conditions $\partial_x u(0,t) = \partial_x u(1,t) = 0$, and initial condition $u(x,0) = 1/(1 + \exp(30x - 10))$. We encode the Neumann conditions via ghost-point extensions $y_0 := y_1$, $y_{N+1} := y_N$, yielding the modified Laplacian

$$[L]_{ij} = \frac{1}{(\Delta x)^2} \cdot \begin{cases} -1 & \text{if } i = j \in \{1, N\} \\ -2 & \text{if } 1 < i = j < N \\ 1 & \text{if } i = j \pm 1 \\ 0 & \text{otherwise,} \end{cases} \tag{32}$$

and the semi-linear ODE $\dot{y}(t) = DLy(t) + R(y(t))$ with element-wise reaction $[R(y)]_i = y_i(1 - y_i)$.

**Further details.** The fixed-step experiment of Section 4.3 uses $N = 100$ with an IWP(2) prior and step grid $\Delta t \in \{10^{-k/2}\}_{k=0}^{7}$.

### A.3.5 2D Fisher-KPP

The 2D Fisher-KPP equation is the two-dimensional analogue of Section A.3.4:

$$\partial_t u(x,y,t) = D\nabla^2 u + u(1-u), \qquad (x,y)\in\Omega, \qquad t\in[0,1], \tag{33}$$

with diffusion $D = 1$, zero-Neumann boundary conditions on $\partial\Omega$, and initial condition $u(x,y,0) = \exp(-30(x^2+y^2))$ seeding an invasion front near the origin. We use a cell-centered $N\times N$ grid with $N = 64$ and centers $(x_i, y_j) = ((i-1/2)/N, (j-1/2)/N)$, stacked column-major as $k = (j-1)N + i$, giving dimension $d = N^2 = 4096$. The 5-point Neumann Laplacian acts as

$$[Ly]_k = \frac{1}{(\Delta x)^2}\sum_{l\in\mathcal{N}(k)} (y_l - y_k), \tag{34}$$

where $\mathcal{N}(k)$ is the set of in-grid 4-neighbors of cell $k$ (4 interior, 3 edge, 2 corner) — equivalent to ghost-point extension of boundary values.

**Further details.** The adaptive work-precision experiment of Section 4.4 uses an IWP(3) prior and the WP tolerance grid.

### A.3.6 2D FitzHugh-Nagumo

The 2D FitzHugh-Nagumo model [19] is a two-species excitable-media PDE on $\Omega = (0,1)^2$:

$$\begin{aligned}\partial_t u(x,y,t) &= a\nabla^2 u + u - u^3 - v + k,\\ \partial_t v(x,y,t) &= (b\nabla^2 v + u - v)/\tau, \qquad (x,y)\in\Omega, \qquad t\in[0,20],\end{aligned} \tag{35}$$

with parameters $a = 2.8\times10^{-4}$, $b = 5\times10^{-3}$, $k = -5\times10^{-3}$, $\tau = 0.1$, zero-Neumann boundary conditions on $\partial\Omega$, and a uniformly random initial condition on $[0,1]^{2N^2}$ (deterministic seed). We discretize both species on a cell-centered $N\times N$ grid with spacing $\Delta x = 1/N$, using the same 5-point Neumann Laplacian as in Section A.3.5. The resulting state $[u_1,...,u_{N^2}, v_1,...,v_{N^2}]$ has dimension $d = 2N^2$.

**Further details.** We evaluate three scaling scenarios, all using an IWP(3) prior without smoothing and adaptive step-size selection at $\tau_{\text{rel}} = 10^{-1}$, $\tau_{\text{abs}} = 10^{-3}$.

- *Domain growth* (weak scaling): $\Delta x = 1/64$ is fixed and $N \in \{8,12,16,24,32,48,64,96,128,192,256\}$ is varied; per-cell stiffness is constant while $d$ grows.
- *Stiffening*: $N = 8$ $(d = 128)$ is fixed and the physical domain size $L \in \{1, 0.5, 0.25, 0.125, 0.0625, 0.03\}$ is varied so that $\Delta x = L/N$ changes; the stiffness proxy $(b/\tau)\Delta x^{-2} = 0.05\Delta x^{-2}$ spans approximately 3–3200.
- *Refinement* (strong scaling): $L = 1$ is fixed and $N \in \{8,12,16,24,32,48,64,96,128\}$ is varied, so that both $d$ and stiffness grow simultaneously as the grid is refined.

Additionally, we visualize the `MatfreeEK1` solution on a $250\times250$ grid ($d = 125,000$) at $\tau_{\text{rel}} = 10^{-2}$, $\tau_{\text{abs}} = 10^{-3}$.